\documentclass{amsart}
\usepackage[dvips]{graphicx}
\usepackage{amscd}
\usepackage{amsmath}
\usepackage{amsfonts}
\usepackage{amssymb}

\usepackage[all]{xy}

\newtheorem{theorem}{Theorem}

\newtheorem{lemma}[theorem]{Lemma}
\newtheorem{proposition}[theorem]{Proposition}
\theoremstyle{definition}

\newtheorem{example}[theorem]{Example}

\theoremstyle{remark}
\newtheorem*{acknowledgements}{Acknowledgements}

\def\openone%{\hbox{\upshape \small1\kern-3.3pt\normalsize1}}
{\mathchoice
{\hbox{\upshape \small1\kern-3.3pt\normalsize1}}
{\hbox{\upshape \small1\kern-3.3pt\normalsize1}}
{\hbox{\upshape \tiny1\kern-2.3pt\SMALL1}}
{\hbox{\upshape \Tiny1\kern-2pt\tiny1}}}
\makeatletter
\newbox\ipbox
\newcommand{\ip}[2]{\left\langle #1\mathrel{\mathchoice
{\setbox\ipbox=\hbox{$\displaystyle \left\langle\mathstrut #1#2\right\rangle$}
\vrule height\ht\ipbox width0.25pt depth\dp\ipbox}
{\setbox\ipbox=\hbox{$\textstyle \left\langle\mathstrut #1#2\right\rangle$}
\vrule height\ht\ipbox width0.25pt depth\dp\ipbox}
{\setbox\ipbox=\hbox{$\scriptstyle \left\langle\mathstrut #1#2\right\rangle$}
\vrule height\ht\ipbox width0.25pt depth\dp\ipbox}
{\setbox\ipbox=\hbox{$\scriptscriptstyle \left\langle\mathstrut
#1#2\right\rangl
e$}
\vrule height\ht\ipbox width0.25pt depth\dp\ipbox}
} #2\right\rangle}
\newcommand{\diracb}[1]{\left\langle #1\mathrel{\mathchoice
{\setbox\ipbox=\hbox{$\displaystyle \left\langle\mathstrut #1\right.$}
\vrule height\ht\ipbox width0.25pt depth\dp\ipbox}
{\setbox\ipbox=\hbox{$\textstyle \left\langle\mathstrut #1\right.$}
\vrule height\ht\ipbox width0.25pt depth\dp\ipbox}
{\setbox\ipbox=\hbox{$\scriptstyle \left\langle\mathstrut #1\right.$}
\vrule height\ht\ipbox width0.25pt depth\dp\ipbox}
{\setbox\ipbox=\hbox{$\scriptscriptstyle \left\langle\mathstrut #1\right.$}
\vrule height\ht\ipbox width0.25pt depth\dp\ipbox}
}\right. }
\newcommand{\dirack}[1]{\left. \mathrel{\mathchoice
{\setbox\ipbox=\hbox{$\displaystyle \left.\mathstrut #1\right\rangle$}
\vrule height\ht\ipbox width0.25pt depth\dp\ipbox}
{\setbox\ipbox=\hbox{$\textstyle \left.\mathstrut #1\right\rangle$}
\vrule height\ht\ipbox width0.25pt depth\dp\ipbox}
{\setbox\ipbox=\hbox{$\scriptstyle \left.\mathstrut #1\right\rangle$}
\vrule height\ht\ipbox width0.25pt depth\dp\ipbox}
{\setbox\ipbox=\hbox{$\scriptscriptstyle \left.\mathstrut #1\right\rangle$}
\vrule height\ht\ipbox width0.25pt depth\dp\ipbox}
} #1\right\rangle}
\let\subsubsubsectionname\@empty
\newcounter{subsubsubsection}[subsubsection]

\def\l@subsubsubsection{\@tocline{4}{0pt}{1pc}{9pc}{}}

\def\subsubsubsection{\@startsection{subsubsubsection}{4}%
  \z@{.5\linespacing\@plus.7\linespacing}{-.5em}%
  {\normalfont\itshape}}
\AtBeginDocument{%
  \@for\@tempa:=-1,0,1,2,3,4\do{%
    \@ifundefined{r@tocindent\@tempa}{%
      \@xp\gdef\csname r@tocindent\@tempa\endcsname{0pt}}{}%
  }%
}
\def\@writetocindents{%
  \begingroup
  \@for\@tempa:=-1,0,1,2,3,4\do{%
    \immediate\write\@auxout{%
      \string\newlabel{tocindent\@tempa}{%
        \csname r@tocindent\@tempa\endcsname}}%
  }%
  \endgroup}
\makeatother
\def\LaTeXparent#1{}%
\def\ChildStyles#1{}%

\newcommand{\frogleg}{\ }

\begin{document}

\title[Non-stationarity of isomorphism]{Non-stationarity of
isomorphism between AF
algebras defined by stationary Bratteli diagrams}
\author{Ola~Bratteli}
\address{Mathematics Institute\\
University of Oslo\\
PB 1053 Blindern\\
N-0316 Oslo\\
Norway}
\email{bratteli@math.uio.no}
\author{Palle~E.~T.~ J\o{}rgensen}
\address{Department of Mathematics\\
The University of Iowa\\
Iowa City, IA 52242-1419\\
U.S.A.}
\email{jorgen@math.uiowa.edu}
\author{Ki Hang Kim}
\address{Mathematics Research Group\\
Alabama State University\\
Montgomery AL~36101--0271\\
U.S.A.}
\email{kkim@asu.alasu.edu}
\author{Fred Roush}
\address{ibid.}
\email{froush@asu.alasu.edu}

\begin{abstract}
We first study situations where the stable AF-algebras defined 
by two square primitive nonsingular incidence matrices 
with nonnegative integer matrix elements are
isomorphic even though no powers of the associated 
automorphisms of the
corresponding dimension groups are isomorphic.
More generally we consider neccessary and
sufficient conditions for two such matrices to
determine isomorphic dimension groups.
We give several examples.
\end{abstract}

\maketitle
%\tableofcontents
%\listoffigures
%\listoftables

\renewcommand{\frogleg}{\\}

This paper was motivated by attempts in \cite{BJO98}
to classify certain AF algebras defined by constant incidence
matrices. The specific incidence matrices considered in \cite{BJO98}
are of the form (18) below, and we shall see there that the first
problem referred to in the abstract is most interesting
for those matrices. The second problem referred to
in the abstract is significant not only for AF algebras but also
for e.g.

   -classification of substitution minimal systems up to strong
    orbit equivalence, \cite{GPS95}, \cite{For97}, \cite{DHS}.

   -homeomorphism classification of domains of certain inverse limit
    hyperbolic systems, \cite{BD95}, \cite{SV98}. The latter paper, 
    which was written independently of this paper, and which was
    pointed out to us by the referee, makes contributions in the same
    direction as our paper. Our $C^*$-equivalence of matrices correspond
    to weak equivalence of (the transposed) matrices in that paper.
    Theorem 2.3 (which is \cite[Corollary 3.5]{BD95}) and 2.4
    in \cite{SV98} corresponds more or less to our Theorem 10.   
    Their Theorem 3.2 is similar to
    our Theorem 6. While the latter part of their paper is
    focused on a class of incidence matrices  arising from
    periodic kneading sequences, our focus here and in \cite{BJO98} 
    is on matrices of the form (18) below which arises in the
    representation theory of Cuntz algebras.

   -cohomology of subshifts of finite type, \cite{BH96}, \cite{Po89}.

In a forthcoming paper we will show  that the isomorphism 
problem for stationary AF algebras is decidable. It is already
known that shift equivalence is decidable in this setting,
\cite{KR79},\cite{KR88}. This is interesting in
view of the fact that isomorphism between two AF algebras 
is known not to be decidable in general, i.e. there is 
no recursive algorithm to decide if two given effective presentations of Bratteli diagrams yield equivalent diagrams in the general (non-stationary) case, 
see \cite{MP98}.

First we will survey some 
terminology and basic facts in the fields of operator
algebras and symbolic dynamics.
Recall from \cite{Bra72} that a $C^{\ast}$-algebra
$\mathfrak{A}$ is called AF (approximately finite dimensional) if it
is the closure of the union an increasing
sequence $\mathfrak{A}_n$ of finite dimensional subalgebras. It is
known from [Bra72,
Theorem~2.7] that two AF algebras
$\mathfrak{A}=\overline{\bigcup_n\mathfrak{A}_n}$,
$\mathfrak{B}=\overline{\bigcup_n\mathfrak{B}_n}$ are isomorphic if
and only if there
are increasing sequences $k_i$, $l_i$ of natural numbers and
injections
$\alpha_i: \mathfrak{A}_{k_i}\to \mathfrak{B}_{\ell_i}$,
$\beta_i:\mathfrak{B}_{\ell_i}\to \mathfrak{A}_{k_{i+1}}$ such that
the following
diagram commutes
\begin{equation}
\vcenter{\xymatrix @-.7pc {
\mathfrak{A}_{k_1} \ar @{_{(}->}[dd] \ar[dr]^{\alpha_1} & \\
   & \mathfrak{B}_{\ell_1} \ar[dl]_{\beta_1} \ar @{_{(}->}[dd] \\
\mathfrak{A}_{k_2} \ar @{_{(}->}[dd] \ar[dr]^{\alpha_2} & \\
   & \mathfrak{B}_{\ell_2} \ar[dl]_{\beta_2} \ar @{_{(}->}[dd] \\
\mathfrak{A}_{k_3} \ar[dr]^{\alpha_3} & \\
   & \mathfrak{B}_{\ell_3}
}}
\end{equation}

\noindent
(The if-part is trivial). This is easily translated into the fact
that there exists a complete isomorphism invariant for AF algebras
$\mathfrak{A}$, namely the dimension group, \cite{Ell76}. 
In the case that
$\mathfrak{A}$ has a unit this is the triple $\left(  K_{0}\left(
\mathfrak{A}\right) ,K_{0}\left(  \mathfrak{A}\right)  _{+},\left[
\openone\right]  \right)  $ where $K_{0}\left(  \mathfrak{A}\right)
$ is an abelian group, $K_{0}\left(
\mathfrak{A}\right)  _{+}$ are the positive elements in
$K_{0}\left(
\mathfrak{A}\right)  $ relative to an order making $K_{0}\left(
\mathfrak {A}\right)  $ into a Riesz ordered group without
perforation, and $\left[
\openone\right]  $ is the class of the identity in $K_{0}\left(
\mathfrak {A}\right)  $ (if $\mathfrak{A}$ is nonunital, replace
$\left[  \openone
\right]  $ by the hereditary subset $\{\left[  p\right]  \mid p$
projection in $\mathfrak{A}\}$ of
$K_{0}\left(  \mathfrak{A}\right) _{+}$). It is also costumary
to apply the term $dimension$ $group$ for just the couple
$\left(  K_{0}\left(
\mathfrak{A}\right) ,K_{0}\left(  \mathfrak{A}\right)  _{+}\right)$,
and this will be done in the sequel. (In dynamical systems theory
this term is used slightly differently; see the comments 
prior to Lemma 1, below.)
See \cite{Eff81}, as well as \cite{LM95},
\cite{Wal92}, \cite{Tor91}, \cite{BMT87}, for
details on this and the following. Let us now specialize to
the case that
$\mathfrak{A}$ is given by a constant
$N\times N$ incidence matrix
$J$ (with nonnegative integer entries) which is primitive, i.e.,
$J^{n}$ has only positive entries for some $n\in\Bbb{N}$,
\cite[Theorem~4.5.8]{LM95}. Then $\mathfrak{A}$ is simple with a
unique trace state $\tau$. In the case that
$K_0(\mathfrak{A})\cong\mathbb{Z}^N$, or, equivalently, when $J$ is
unimodular, this class of AF algebras (or rather dimension groups)
has been characterized intrinsically in [Han81, Theorems~3.3 and
4.1]. We do not assume unimodularity in the sequel.

\vskip0.3cm

In general when $J$ is an
$N\times N$ matrix with nonnegative entries, the dimension group is
the inductive limit
\begin{equation}
\mathbb{Z}^{N}\overset{J}{\longrightarrow}\mathbb{Z}^{N}\overset
{J}{\longrightarrow}\mathbb{Z}^{N}\overset{J}{\longrightarrow}\cdots
\end{equation}
with order generated by the order defined by
\begin{equation}
\left(  m_{1},\dots,m_{n}\right)  \geq0\Longleftrightarrow m_{i}
\geq0\text{\qquad on }\mathbb{Z}^{N}.
\end{equation}
This group can be computed explicitly as a subgroup of
$\mathbb{Q}^{N}$ as follows when $\det\left(  J\right)  \neq0$ (as
it will be in our examples): Put
\begin{equation}
G_{m}=J^{-m}\left(  \mathbb{Z}^{N}\right)  ,\qquad m=0,1,\dots,
\end{equation}
and equip $G_{m}$ with the order
\begin{equation}
G_{m}^{+}=J^{-m}\left(  \left(  \mathbb{Z}^{N}\right)  ^{+}\right)  .
\end{equation}
Then
\begin{equation}
G_{0}\subset G_{1}\subset G_{2}\subset\cdots,
\end{equation}
and
\begin{equation}
K_{0}\left(  \mathfrak{A}_{L}\right)  =\bigcup_{m}G_{m},
\end{equation}
is a subgroup of $\mathbb{Q}^{N}$ (containing $\mathbb{Z}^{N}$),
with order defined by
\begin{equation}
g\geq0\text{\qquad if }g\geq0\text{ in some }G_{m}.
\end{equation}
The action of the trace state $\tau$ on $K_{0}\left(
\mathfrak{A}\right)
$ may be computed as follows: If $\lambda$ is the Frobenius
eigenvalue of
$J$, and
$\alpha=\left(  \alpha_{1},\dots,\alpha_{N}\right) $ is a
corresponding eigenvector
in the sense
\begin{equation}
\alpha J=\lambda\alpha
\end{equation}
(i.e., $J^{t}\alpha^{t}=\lambda\alpha^{t}$, see
\cite[pp.~33--37]{Eff81}), then if $\alpha$ is suitably normalized
(by multiplying with a positive factor), the trace applied to an
element $g$ at the $m$'th stage of
\begin{equation}
\underset{1}{\mathbb{Z}^{N}}
\longrightarrow\underset{2}{\mathbb{Z}^{N}
}\longrightarrow\underset{3}{\mathbb{Z}^{N}}\longrightarrow\dots
\longrightarrow\overset{
\begin{array}[c]{c}
g\\
\makebox[0pt]{\hss$\displaystyle\cap$\hss}\makebox[0pt] {\hss\rule
[-0.15pt]{0.225pt}{6pt}\hss}
\end{array}
}{\underset{m}{\mathbb{Z}^{N}}}\longrightarrow\cdots
\end{equation}
is
\begin{equation}
\tau\left(  g\right)  =\lambda^{-m+1}
\langle \alpha | g\rangle,
\end{equation}
where $\langle \cdot | \cdot\rangle$
 here denotes the usual inner product in
$\mathbb{R}^{N}$, i.e., $\langle \alpha |
g\rangle=\sum_{i=1}^{N}\alpha_{i}g_{i}$.
Taking
$\alpha$ as the Frobenius eigenvector makes the ansatz well
defined: if $g\in G_{m}\subset G_{m+1}$, then
\begin{equation}
\lambda^{-m+1}\langle \alpha | g\rangle=\lambda^{-\left(  m+1\right)
+1}\langle \alpha | Jg\rangle .
\end{equation}
Thus $\tau$ is an additive character on $K_{0}\left(
\mathfrak{A}\right)$, and up to normalization the unique positive
such. If
$\mathfrak{A}$ is unital we
may normalize $\alpha$ by requiring
$\langle \alpha | [\openone]_0\rangle=1$, and it can then be shown
that the range of
the trace on projections is
$\tau\left( K_{0}\left( \mathfrak{A}\right) \right)  \cap\left[
0,1\right]  $.

When $K_{0}\left( \mathfrak{A}\right) $ is given concretely in
$\mathbb{Q}^{N}$ as above, the trace can be computed as
\begin{equation}
\tau\left(  g\right)  = \langle \alpha | g\rangle,
\end{equation}
where $g\in m$'th term $\mathbb{Z}^{N}$ is identified with its
image $J^{-m+1}g$ in $\mathbb{Q}^{N}$; and the positive cone in
$K_{0}\left( \mathfrak{A}\right)  \subset\mathbb{Q}^{N}$ identifies
with those $g$ such that $\tau\left( g\right) >0$, or $g=0$.

If one forgets about $\left[\openone\right]$, or the hereditary
subset of $K_0(\mathfrak{A})_+$, then $(K_0(\mathfrak{A}),
K_0(\mathfrak{A})_+)$ is a complete
invariant for stable isomorphism, i.e. isomorphism of
$\mathfrak{A}\otimes\mathcal{K}(L^2(\mathbb{Z}))$. 
In the rest of this paper we will only consider
stable isomorphism, and not study the position of 
$\left[\openone\right]$
inside $(K_0(\mathfrak{A}), K_0(\mathfrak{A})_+)$. 
Let us just mention that in the applications in \cite{BJO98},
the element $\left[\openone\right]$ is represented by
$(1,0,0,\dots,0)^T$ in the concrete representation (7),
and to take care of $\left[\openone\right]$ one has to
assume that the matrices $A_i, B_i$ in (15) below
preserves the class of $\left[\openone\right]$. In Theorems
6 and 7 below this amounts to the added condition that
$$
A(1)J^{n_0}(1,0,0,\dots,0)^T = K^{m_0}(1,0,0,\dots,0)^T
$$
for some non-negative integers $n_0, m_0$. For more details on isomorphisms as opposed to stable isomorphisms in this setting, see
\cite{BJO98}.

If $\mathfrak{A},\mathfrak{B}$ are AF algebras defined by constant
(necessarily square) incidence matrices
$J, K$, it
follows from \cite[Theorem 2.7]{Bra72} that $\mathfrak{A}$ and
$\mathfrak{B}$ are
stably isomorphic if and only if there exist natural numbers
$n_{1},n_{2},n_{3},\dots$, $m_{1},m_{2},m_{3},\dots$, and matrices
$A_{1},A_{2},\dots$, $B_{1},B_{2},\dots$ with nonnegative integer
matrix elements such that the following diagram commutes:
\begin{equation}
\vcenter{
\xymatrix @-.7pc {
\bullet \ar[dd]_{J^{n_1}} \ar[drr]^{A_1} & & \\
&   & \bullet \ar[dll]_{B_1} \ar[dd]^{K^{m_1}} \\
\bullet \ar[dd]_{J^{n_2}} \ar[drr]^{A_2} & \\
&   & \bullet \ar[dll]_{B_2} \ar[dd]^{K^{m_2}} \\
\bullet \ar[dd]_{J^{n_3}} \ar[drr]^{A_3} & \\
&   & \bullet \ar[dll]_{B_3} \ar[dd]^{K^{m_3}} \\
\bullet \ar[dd] \ar[drr]^{A_4} & \\
&   & \bullet \ar[dll] \ar[dd] \\
\vdots & \\
&   & \vdots
}}
\end{equation}
This means that
\begin{equation}
\begin{array}{l}
J^{n_{k}}   =B_{k}A_{k},\\
K^{m_{k}}   =A_{k+1}B_{k}\end{array}
\end{equation}
for $k=1,2,\dots$. The first aim of this paper is to show that the
sequences $A$, $B$, $n$, $m$ cannot in general be taken to be
constant when they exist. 

\vskip0.3cm

In the covariant version of this isomorphism problem, it is known
from a theorem of Krieger that the sequences can be taken to be
constant. Let $G\left(  J\right)  $ be the dimension group
associated to $J$, and $\left( \sigma_{J}\right)  _{\ast}$ the shift
automorphism of $G\left( J\right)  $ determined by $J$,
\cite[pp.~36--37]{Eff81}. Let now $\mathfrak{A}$ be the stable
AF-algebra associated to $G\left(  J\right)  $, and $\sigma_{J}$ an
automorphism of
$\mathfrak{A}$ such that the corresponding automorphism of $G\left(
J\right)  $ is $\left(  \sigma_{J}\right)  _{\ast}$. Then Krieger's
theorem \cite{Kri80} says that $\left(  G\left(  J\right)  ,\left(
\sigma_{J}\right) _{\ast}\right)  $ is isomorphic to $\left(
G\left(  K\right) ,\left(  \sigma_{K}\right)  _{\ast}\right)  $ if
and only if there is a $k\in\mathbb{N}$ and nonnegative rectangular
matrices $A$, $B$ such that
\begin{equation}
\begin{array}{l}AJ  =KA,\\
BK   =JB,\\
BA   =J^{k},\\
AB   =K^{k}. \end{array}
\end{equation}
If also $N>1$, it was proved recently in \cite{BrKi98} that this is
also equivalent to outer conjugacy of $\sigma_J$ and $\sigma_K$.
(This was proved in \cite{EvKi97} in the case that
$(K_0(\mathfrak{A}),K_0(\mathfrak{A})_+)$ has no infinitesimal
elements, i.e. this
ordered group is totally ordered.)

So, in dynamical system language, the problem is:
Given primitive square matrices
$J,K$ over the nonnegative integers such that there exist
sequences $m(i),n(i)\in \mathbb{N}$ and matrices $A(i),B(i)$ over
the non-negative
integers for $i\in \mathbb{N}$ with
\begin{equation}
\begin{array}{l} J^{n(i)}=B(i)A(i) , \\
  K^{m(i)}=A(i+1)B(i) ,\end{array}
\end{equation}
are some positive powers of $J,K$ elementary shift equivalent?

We will show in Proposition 2 and Proposition 5 that the answer of this
question is no in general. The difference between these two propositions
is that the matrices in Proposition 5 has the special form (18) below.
Note that the pair of matrices in Proposition 5 are not 
unimodular, and
\[
G_m\subsetneqq G_{m+1},
\]
where $G_m$ is defined by (4). Moreover, the two matrices
have the form (18) which was
the one required for the \cite{BJO98}
analysis. In general, we have
$(G_{m+1}:G_m)=|\det J|$, and for
the $J,K$ pair in Proposition~5,
$G_{m+1}\diagup G_m$ is $\mathbb{Z}\diagup 32\mathbb{Z}$ for $J$,
and
$\mathbb{Z}\diagup 16\mathbb{Z}$ for $K$.

\vskip0.3cm

A few words about terminology: In the theory of symbolic 
dynamics the term "dimension group" is used slightly differently
from the usage in $C^*$-theory introduced above, namely
for the abelian group $G\left(J\right)$ without order structure, and it
is used even in the wider context of non-positive matrices, 
\cite[Definition 7.5.1]{LM95}, \cite{BMT87}. The shift automorphism 
$\left(  \sigma_{J}\right)  _{\ast}$ is called the $dimension$ 
$group$ $automorphism$ in this context. If $J$ is non-negative,
the positive part $G\left( J\right)  _{+}$ is called the
$dimension$ $semigroup$, and the triple $\left(  G\left(J\right), 
 G\left( J\right)  _{+},  \left(  \sigma_{J}\right)  _{\ast} \right)  $
is called the $dimension$ $triple$. In the rest of this paragraph,
let the term "matrix" mean "matrix over the non-negative integers".
Then two square matrices $J,K$ are $elementary$ $shift$ 
$equivalent$ if there exist matrices $A,B$ such that $J=BA$
and $K=AB$. We say that $J,K$ are $shift$ $equivalent$ of $lag$ $k$
if (16) holds for some matrices $A,B$, and they are $shift$
$equivalent$ if they are shift equivalent of some lag $k$ in
$\mathbb{N}$. Thus shift equivalence of lag 1 is the same
as elementary shift equivalence, \cite[Proposition 7.3.2]{LM95}.
The matrices $J,K$ are $strongly$ $shift$ $equivalent$ if they can
be connected by a finite chain of elementary shift equivalent 
matrices. (Note that elementary shift equivalence is not an 
equivalence relation; it is not transitive.) Strong shift equivalence
trivially implies shift equivalence, \cite[Theorem 7.3.3]{LM95},
but the converse is the long standing Williams conjecture, and it is not
true, \cite{KR92}, even when $J$ and $K$ are irreducible, 
\cite{KiRo98}. Since the stable AF algebras defined by $J,K$ are isomorphic if and only if (17) holds for some sequences $A,B,n,m$,
the equivalence defined by (17), or (14)-(15), could be 
termed $C^*$-$equivalence$.

\vskip0.3cm

In conclusion, we consider in this paper the following notions 
of equivalence of two non-negative matrices $J,K$, where each
notion is strictly stronger than the next one. The notions 2, 3, 4 and 5
are equivalence relations, and 2 is the equivalence relation
generated by 1.

\vskip0.3cm

1. Elementary shift equivalence = shift equivalence of lag 1.

\vskip0.3cm

2. Strong shift equivalence

\vskip0.3cm

3. Shift equivalence = isomorphism of dimension triples (by 
Krieger«s theorem)

\vskip0.3cm

4. Shift equivalence of some positive powers $J^n$ and $K^m$ =
elementary shift equivalence of some powers

\vskip0.3cm

5. $C^*$-equivalence = isomorphism of (ordered) dimension groups = 
isomorphism of associated stable AF algebras

\vskip0.3cm

(It is well known that 2, 3 and 5 are equivalence relations,
and for 4 the argument is as follows: If $J,K$ and $K,L$ are
related as in 4, then there exist positive powers 
$m,n,k,j$ such that each of the pairs $J^n,K^m$ and $K^k,L^j$
are elementary shift equivalent. But then each of the pairs
$J^{nk},K^{mk}$ and $K^{mk},L^{mj}$ are elementary shift equivalent,
and thus shift equivalent. Since shift equivalence is transitive,
it follows that $J^{nk},L^{mj}$ are shift equivalent.)

\vskip0.3cm

Recall that elementary shift equivalence of non-singular
square matrices implies conjugacy in the usual matrix sense 
over $\mathbb{Q}$, even in the absense of positivity.
This is even true for shift equivalence, by an argument in the
next paragraph. Note also that if $J,K$ are primitive matrices,
then shift equivalence in the sense of (16) is equivalent to
the corresponding relation where $A,B$ are only assumed 
to be matrices over $\mathbb{Z}$ (not assuming non-negativity), 
\cite{PW77}, \cite[Theorem 2.1]{KR79}, \cite[Theorem 7.3.6]{LM95}.

\vskip0.3cm

The problem addressed in this paper arose in \cite{BJO98}
for incidence matrices $J$ and $K$ of the form  
\begin{equation}
J=\left( \begin{array}{cccccc}
m_1 & 1 & 0 & \cdots & 0 & 0\\
m_2 & 0 & 1  & \cdots & 0 & 0\\
\vdots & \vdots & \vdots & \vdots & \vdots & \vdots \\
m_{N-1} & 0 & 0 & \cdots & 0 & 1 \\
m_N & 0 & 0 & \cdots & 0 & 0\end{array}\right)
\end{equation}
where $m_i$ are non-negative integers, $m_N\not=0$ and the greatest common divisor of the set of $k$
such that $m_k\not=0$ is 1.  The last conditions
ensure that these matrices are nonsingular and primitive,
and whenever we refer to a "matrix of type (18)" we assume that
these additional conditions are satisfied. At the outset,
it was not even clear if there were different such matrices with
isomorphic $C^*$-algebras. The paper \cite{BJO98} does, 
however, contain
examples of pairs $J,K$ of distinct matrices of type (18)
which are $C^*$-equivalent.
It is ironic that while the present paper started in a quest
for pairs of the type (18) defining stably isomorphic
AF algebras, i.e., pairs of matrices satisfying the 
condition (17), but not with any constant sequences $A, B, n, m$, 
it is not so easy to find
an example of a pair of distinct such matrices satisfying
the condition {\em with} constant $A, B, n, m$. We give such an
example between (32) and (33), and another example which is
close in that $J^{12}$ is conjugate to $K^{12}$ in (31).
Note in this connection that two non-negative square 
matrices $J,K$ of the form (18) which are shift equivalent 
are identical by the following reasoning: They are nonsingular
and primitive, and since they are $C^*$-equivalent they have the
same size $N\times N$, where $N$ is the rank of the associated
dimension group. The two last relations in (16) then imply 
that $A$ and $B$ are nonsingular $N\times N$ matrices. Hence
any of the first two relations in (16) implies that $J$ and $K$ are
conjugate over $\mathbb{Q}$, and thus they 
have the same characteristic polynomial.
But as explained in the beginning of Example 9, below, the
characteristic polynomial uniquely determines these matrices,
and hence $J=K$. Thus for matrices $J,K$ of type (18),
shift equivalence is the same as equality, $C^*$-equivalence
is $not$ the same as equality (by Proposition 5 and Example 9),
while the situation for the intermediate equivalence 4. is
also that shift equivalence of some powers do not imply equality,
see the penultimate example in Example 9 , or [BJO98, Section~7]. 
(The matrices (31) are 
$C^*$-equivalent with conjugate twelfth powers, but 
apparently do not satisfy 4.)

\begin{lemma}
Suppose that two given non-negative square $d\times d$ 
matrices $J,K$ are primitive,
nonsingular and equal at their largest eigenvalue,
that is, they have the same Perron eigenvalue and
 row and column eigenvectors.
Then there exists a positive integer $c$ such that $J^{cn}K^{-n}$ is
nonnegative for all $n>0$.
\end{lemma}

\begin{proof}
The  parts at the largest
eigenvalue will multiply, be positive,
and will swamp all the others, since they grow at
an exponential rate corresponding to this eigenvalue.
More precisely, we can conjugate and then write
$J=J_1+J_2,K=K_1+K_2$
corresponding to the eigenspace for the maximal
eigenvalue $\lambda_1$,
and the eigenspaces for all other eigenvalues.
Let the maximum of the absolute values of
those eigenvalues be $\lambda_2$
and the maximum absolute value for an eigenvalue of
their inverses be say $\lambda_3$.  Then
\begin{equation}
J^{ck}K^{-k}=J_1^{ck}K_1^{-k} + J_2^{ck}K_2^{-k}.
\end{equation}
After we conjugate back, the entries in the first summand
contribute entries proportional to the fixed row
 and column eigenvectors, which
are at least $C_1\lambda_1^{ck-k}$.
The entries in the 2nd part are at most
$C_2 \lambda_2^{ck}\lambda_3^k$.
 Choose $c$ large enough that
 $\lambda_1^{c-1}>\lambda_2^c \lambda_3$
and the inequality will eventually hold.
\end{proof}

Let us now consider the following two matrices
\begin{eqnarray}
&&J=\left( \begin{array}{ccccc}
1 &1 &0 &0 &0 \\
     0 &1 &1 &0 &0 \\
     0 &0 &1 &1 &0 \\
     0 &0 &0 &1 &1 \\
     1 &0 &0 &0 &1\end{array} \right)=1+P\;, \\ [.5ex]
&&K=\left( \begin{array}{ccccc}
0 &1 &1 &0 &0 \\
0 &0 &1 &1 &0 \\
1 &0 &0 &0 &1 \\
1 &1 &0 &0 &0 \\
0 &0 &0 &1 &1\end{array} \right) =QJ\;,
\end{eqnarray}
where $P$ and $Q$ denote the
 matrices of the cyclic permutations which
send 1,2,3,4,5 respectively to 2,3,4,5,1
and 2,3,5,1,4.

\begin{proposition}
The matrices $J, K$ are $C^*$-equivalent, but no two powers
are conjugate over the rational numbers.
\end{proposition}

\begin{proof}
The $C^*$-equivalence is the
existence of an infinite sequence
$A(i),B(i)$ of nonnegative matrices,
and suitable powers, such that
we have
$J^{n(i)}=B(i)A(i)$, $\,K^{m(i)}=A(i+1)B(i)$.
Let $A(1)=I$. We can solve recursively
\begin{equation}
\begin{array}{l}B(1)=J^{n(1)} \\
     A(2)=K^{m(2)}J^{-n(1)}\\
     B(2)=J^{n(2)+n(1)}K^{-m(2)}\\
     A(3)=K^{m(3)+m(2)}J^{-n(2)-n(1)} \end{array}
\end{equation}
and so on.  It follows that if for any $c$
there exists $d$ such that
$J^dK^{-c}, K^dJ^{-c}$ are nonnegative integer matrices, the
results hold.  Nonnegativity follows from Lemma 1:
the matrices are primitive, their Perron
eigenvalue is 2 and their left and right Perron eigenvectors are
$(1,1,1,1,1)$ and
$(1,1,1,1,1)^T$, respectively.
To
 verify integrality, we compute
the determinants as both 2; then the row spaces of
$J^n,K^n$ each lie within the space of vectors $w$
whose product
with the column vector $(1,1,\ldots ,1)^t$ is
 a positive integer multiple of $2^n$.  But since
$2^n$ is their determinant, this is their exact
row spaces, and the same holds
for column spaces.  It follows that each matrix
is the product of the other
matrix and a unimodular matrix, since the rows
 of the powers of each lie
in the row spaces of the powers of the other.

To show powers of $J,K$ can never be conjugate
over the rational numbers
(which the last equation implies)  we compute that their
characteristic polynomials are
$$
\det(t\openone\! -\! J)=(t-2)(t^4\!-3t^3\!+4t^2\!-2t+1), \quad
\det(t\openone\! -\! K)=(t-2)(t^4\!+t^3\!+1)
$$
and the degree 4 factors are irreducible in $\mathbb{Q}[t]$.
  We restrict to the eigenspaces associated with
 the degree 4 factors.
The discriminants of the degree 4 factors and of
their algebraic number fields are 229, 125. These
are relatively prime, and the root field of $J$ is
cyclotomic (5th roots of unity
from its circulant form), its only nontrivial proper subfield
is quadratic and can be determined
also to have discriminant a multiple of 5.
 Therefore the only intersection of
the fields is the rational numbers \cite{W63}.  Suppose
we have powers which are conjugate
over the rational numbers, so have the same
 eigenvalues.  These powers of eigenvalues
all are in the intersection field, the rational
numbers, and their product is 1 by the
determinant.  So powers of the eigenvalues
are rational units, $\pm 1$ and
all eigenvalues of both matrices (other than
 the Perron eigenvalues) must
be roots of unity.  But this is false.
\end{proof}

We will next show that even within the
class of algebras with incidence matrices of the form (18) there are non-stationary isomorphisms. To this end we need a more general
version of Lemma~1.

\begin{lemma}
Suppose the matrices $J,K$ are primitive,
nonsingular, of equal size.  Let $V_1,V_2$ be
the sums of all column eigenspaces for eigenvalues other
than the PF eigenvalue of $J,K$, respectively, and
let $w_1,w_2$ be positive eigenvectors corresponding
to the PF eigenvalues of $J,K$, respectively.
Suppose there is a
nonnegative integer matrix $A(1)$ such that
 \vskip0.3cm
(1) $A(1)V_1= V_2$ 
 \vskip0.3cm
(2) $A(1)w_1,A(1)^{-1}w_2$ project
nontrivially to positive multiples of $w_2,w_1$ in the
direct sum of eigenspaces.
\vskip0.3cm
Then there exists a positive integer $c$ such 
that $J^{cn}A(1)^{-1}K^{-n}$, $K^{cn}A(1)J^{-n}$ are
nonnegative for all sufficiently large $n$.
\end{lemma}

\begin{proof}
The  parts at the largest
eigenvalue will multiply, be positive,
and will swamp all the others, since they grow at
an exponential rate corresponding to this eigenvalue.
More precisely, we can conjugate and then write
$J=J_1+J_2,K=K_1+K_2$
corresponding to the eigenspace for the maximal
eigenvalue $\lambda_1$,
and the eigenspaces for all other eigenvalues.
Let the maximum of the absolute values of
those eigenvalues be $\lambda_2$
and the maximum absolute value for an eigenvalue of
their inverses be say $\lambda_3$.  Then by (1)
\begin{equation}
K^{ck}A(1)J^{-k}=K_1^{ck}A(1)J_1^{-k}+
K_2^{ck}A(1)J_1^{-k}+K_2^{ck}A(1)J_2^{-k}.
\end{equation}
After we conjugate back, the entries in the first summand
contribute entries proportional to the fixed row
 and column eigenvectors, which
are at least $C_1\lambda_1^{ck-k}$.
The entries in the 2nd part and 3rd part are at most
$C_2 \lambda_2^{ck}\lambda_3^k$.
 Choose $c$ large enough that
 $\lambda_1^{c-1}>\lambda_2^c \lambda_3$
and the inequality will eventually hold.
The same holds true for the $A(1)^{-1}$ term.
\end{proof}

\begin{lemma}
If for all primes $p$ dividing $det(J)det(K)$
the nonsingular matrices $J,K$ are nilpotent modulo $p$
and the conditions of Lemma~3 are true for some
$A(1)$ whose determinant involves only these primes,
then the equations {\em (15)}
hold for $A(1)$ and appropriate sequences
$n(k),m(k)$.  (This includes the case when there are no
such $p$).
\end{lemma}

\begin{proof}
We write out the equations for a general $A(1)$
which is assumed to be a nonnegative unit
\begin{equation}
\begin{array}{l}
     B(1)=J^{n(1)}A(1)^{-1} \\
     A(2)=K^{m(1)}A(1)J^{-n(1)} \\
     B(2)=J^{n(2)+n(1)}A(1)^{-1}K^{-m(1)} \\
     A(2)=K^{m(2)+m(1)}A(1)J^{-n(2)-n(1)}\;. \end{array}
\end{equation}
The nilpotence modulo $p$ guarantees that some powers of
$J,K$ are divisible by $p$, hence any sufficiently
large power are divisible by the determinant of a given
power of the other matrix, so that if the powers increase
sufficiently rapidly and $A(1)$ has determinant dividing
some power of the determinants of $J,K$, these matrices
exist over the integers, and positivity follows from
Lemma 3 (except for the first equation, which follows
by a similar argument or can be checked step by step).
\end{proof}

\begin{proposition}
The matrices below  are $C^*$-equivalent, but no powers are shift
equivalent over the rationals.
\begin{equation}
\left( \begin{array}{cc} 4 &1 \\ 32 &0\end{array}\right),\qquad
\left( \begin{array}{cc} 6 &1 \\ 16 &0\end{array} \right).
\end{equation}
\end{proposition}

\begin{proof}
We check that the condition in Lemma~4 holds 
with $A(1)=I$ and that the matrices
have eigenvalues respectively $8,-4$; $8,-2$.  At the
negative eigenvalues both column eigenvectors are
$(1,-8)^T$ so the identity maps one to the other.
One can alternatively check by the recipe (4) - (13) that the two 
matrices define the same dimension group (see below).
But the values of the two pairs of eigenvalues prevents any 
power of one matrix to be conjugate to a power 
of the other matrix over the rationals.
\end{proof}

Note that in both the examples
$J=\left(\begin{array}{cc}
4 & 1\\ 32 & 0\end{array}\right)$
and $K=\left(
\begin{array}{cc} 6 & 1 \\ 16 & 0\end{array}\right)$ above,
the concrete realization of the dimension group $G$ as a subgroup of
$\mathbb{Q}^2$ is $\mathbb{Z}[1/2]^2$, and
$(x,y)\in\mathbb{Z}[1/2]^2$ is positive and nonzero iff $8x+y>0$.

This is proved as follows: One shows that
\begin{equation}
J=\left(\begin{array}{cc}
4 & 1\\ 32 & 0\end{array}\right)
\end{equation}
has eigenvalues $8, -4$ with left eigenvectors $(8,1),(4,-1)$ and
right eigenvectors $\left(\begin{array}{c}1 \\ 4\end{array}\right),
\left(\begin{array}{c} -1 \\ 8\end{array}\right)$. Hence
\begin{equation}
\begin{array}{l}
J^n =
\left( \begin{array}{ll}
\frac{2}{3}8^n+\frac{1}{3}(-4)^n &
     \frac{1}{12}8^n-\frac{1}{12}(-4)^n \\ [.5ex]
\frac{8}{3}8^n-\frac{8}{3}(-4)^n &
     \frac{1}{3}8^n+\frac{2}{3}(-4)^n \end{array}\right) \\ [2.5ex]
\qquad = \frac{4^{n-1}}{3}\left(\begin{array}{ll}
4(2\cdot 2^n+(-1)^n) & 2^n-(-1)^n \\
32(2^n-(-1)^n) & 4(2^n+2(-1)^n)\end{array}\right) \end{array}
\end{equation}
for all $n\in\mathbb{Z}$. But as $2\equiv-1\pmod{3}$ and hence
$2^n\equiv(-1)^n\pmod{3}$ for $n\in\mathbb{Z}$ one computes that
the matrix elements for $J^n$ are contained in $\mathbb{Z}[1/2]$
for both negative and positive $n$, and it follows that the
dimension group
$G(J)=\bigcup\limits_{n=0}^\infty J^{-n}\mathbb{Z}^2$ is contained
in $\mathbb{Z}[1/2]^2$. But since $J^n$ is $4^{n-1}$ times a matrix
with integer coefficients by the formula above, it follows that
$4^{1-n}(\mathbb{Z}\times\mathbb{Z})\subseteq G(J)$ for
$n=1,2,\ldots$, and hence $\mathbb{Z}[1/2]^2\subseteq G(J)$. Thus
$G(J)=\mathbb{Z}[1/2]^2$. If $(x,y)\in G(J)$ it follows from the
form of the left eigenvector that $(x,y)>0$ if and only if $8x+y>0$.

Correspondingly
\begin{equation}
K=\left( \begin{array}{cc} 6 & 1 \\ 16 & 0\end{array}\right)
\end{equation}
has eigenvalues $8,-2$ with left eigenvectors $(8,1)$ and $(2,-1)$
and right eigenvectors
$\left(\begin{array}{c} 1 \\ 2\end{array}\right)$ and
$\left(\begin{array}{c} -1 \\ 8\end{array}\right)$. Hence
\begin{equation}
K^n=\frac{2^{n-1}}{5}\left(\begin{array}{ll}
2(4\cdot 4^n+(-1)^n) & 4^n-(-1)^n \\
16(4^n-(-1)^n) & 2(4^n+4(-1)^n)\end{array}\right)
\end{equation}
for all $n\in\mathbb{Z}$, and using $4=-1\pmod{5}$ and hence
$4^n=(-1)^n\pmod{5}$ for $n\in\mathbb{N}$, one proceeds as in the
previous case.

\vskip0.3cm

 The criteria in Theorems 6 and 7, below,
reduce the question of isomorphism of
this kind of AF algebra to standard questions in matrix
theory somewhat like those in \cite{KR79}, that is,
existence of a nonnegative
 integer matrix $A(1)$ which maps certain
computable linear spaces associated with $J$ over extension
fields isomorphically to corresponding linear spaces associated with
$K$.  In fairly simple cases computations should be practical,
 and we will show in a forthcoming paper that, as in the 
case of shift equivalence, \cite{KR79}, \cite{KR88}, 
an algorithm exists which will always
decide isomorphism of the algebras. The hypotheses that
$J$ and $K$ are non-singular in these propositions could
be removed with some reformulation and a longer proof.

\begin{theorem}
Let $J, K$ be primitive, nonsingular, nonnegative square
matrices.  Let $V_1,V_2$ be the sums of all column eigenspaces
for eigenvalues other than the PF eigenvalue of $J,K$ respectively,
and let $w_1,w_2$ be positive eigenvectors corresponding to
the PF eigenvalues of $J,K$, respectively.  In order for
the stationary $C^*$-algebras defined by $J,K$ to be isomorphic,
it is necessary that there be a nonnegative integer matrix $A(1)$
such that
\vskip0.3cm
(1) $A(1)(V_1)=V_2$,
 \vskip0.3cm
(2) $A(1)w_1, A(1)w_2$ project nontrivially to positive
multiples of $w_2,w_1$ by the projection from the sum of
eigenspaces to the eigenspace for the maximal eigenvector
for $K,J$, respectively,
\vskip0.3cm
and
 \vskip0.3cm
(3) for the row (left) PF eigenvector
$v_1$ of $J$, $v_1A(1)^{-1}$ is nonnegative. 
\vskip0.3cm
These conditions are necessary for the nonnegativity
of $A(i),i>1, B(i),i>0$, and are sufficient if we
allow replacement of $A(1)$ by some $A(1)J^{s_0}$.
\end{theorem}

\begin{proof}
 If (1) does not hold true, then
$A(1)$ will map vectors from $V_1$ into the maximal column
eigenspace of $K$ nontrivially.
At the largest eigenvalue where
this occurs, these terms will become dominant
in $K^r A(1)J^{-s}$ and give the asymptotic value of the
entire matrix.   This will
make the limit of $K^r A(1)J^{-s}$ as $r,s\rightarrow
\infty$ in any way, a limit of matrices whose column vectors
come from the nonmaximal eigenspace of $K$.  But this
is impossible, since other nonmaximal eigenspaces of a positive
matrix contain no nonnegative vectors (if they did, multiplication
by powers of the matrix would increase them at a rate
which is asymptotically the maximal eigenvalue, which
means that they would have components in the maximal eigenspace).

This proves (1) and given (1) we have
$$K^r A(1) J^{-s}= K_1^r A(1) J_1^{-s}+K_2^r A(1) J_1^{-s}
+K_2^r A(1) J_2^{-s}.$$
If $A(1)w_1$ projects to a negative multiple of $w_2$,
then those terms will be dominant and make the entirety
negative.  If it projects to a zero multiple, then
$A(1)w_1\subset V_2$ which makes equality in (2) impossible.
This proves (2).  Sufficiency of (1),(2) are proved in
Lemma 3.

For the condition (3), it will suffice that $J^{-s}A(1)^{-1}$
is eventually positive, looking at the dominant maximal
eigenspaces.  Conversely, if the vector is negative then
the dominant part is negative, which is impossible.
Suppose it is nonnegative but not positive.  Then the
indicated replacement continues to allow solution
of the other equations with altered exponents, but
it maps the nonnegative vector to a positive one.
\end{proof}

  Define $p$-adic limits of the powers of integer matrices $A$
whose determinants divide $p$ as follows:  modulo each
power $p^m$, there is a unique power $A^{e(m)}$
which is idempotent, since $\{ A^n |n\in \mathbb{Z}^+ \}$
is a cyclic finite semigroup.
These idempotent powers agree to reductions modulo the
lower of the two powers of $p$, by uniqueness.  Therefore
$A^{e(m)}$ have a $p$-adic limit.
 Modulo each power of $p$, any sufficiently large
powers of the original matrix have the same row spaces as
the idempotents, since each is a power of $A$ times the
other.

The result in the next theorem also holds if
$\mathbb{Z}$ is replaced by the ring of algebraic integers of any
algebraic number field, and $p$ by any prime ideal of that ring.
The same proof goes through. Recall here and in Propostion~10
that every algebraic number field contains the unique
subring of algebraic integers, all elements in the field
which satisfy monic polynomial equations.  Primes refers 
in general to prime ideals in this subring (which give rise to
valuations on the field).  Given an element in the
field, we can factor the ideals generated by its numerator 
and denominator 
into prime ideals uniquely, subtract, and hence up to units
write it as a product of positive and negative powers
of primes. These prime ideals will not in general be principal
and hence arise from single elements, but a finite index
subgroup of their multiplicative group are principal ideals; the
finite quotient group is called the (ideal) class group
of the field.

 \begin{theorem}
 In order for a pair of non-singular non-negative 
square matrices $J,K$ to be $C^*$-equivalent,
it is necessary that there exists a non-negative integer matrix
$A(1)$ sending the $p$-adic row space
 of $J_1$ to the $p$-adic row space of
$K_1$ isomorphically, where $J_1,K_1$
are the $p$-adic limits of powers of $J,K$ respectively.
This condition, taken over all primes dividing
the determinants of $J,K,$ is necessary and sufficient that we can
(possibly altering exponents and taking a replacement
of the matrix $A(1)$) make all the matrices
$A(i),B(i)$ in (17), (24) have integral entries.
\end{theorem}

\begin{proof}
  If $A(1)$ gives an isomorphism then for arbitrary large
powers of $J,K$, the matrices
   $K^c A(1) J^{-d}= U,\ J^e A(1)^{-1}K^{-f}=V$
are integral. Modulo any fixed power of $p$, we can arrange
by increasing these powers and altering $U,V$ to
other integral matrices that the powers are in each
case those giving rise to the idempotent limits
$K_1,J_1$.  Therefore modulo each power of $p$,
$$
 K_1 A(1) =UJ_1, J_1=VK A(1).
$$
Hence these matrices $K_1A(1), J_1$ have equal
row spaces.

Conversely, suppose that this condition holds.  Then
we can find $p$-adic $U,V$ satisfying the last equations.
Hence they satisfy them modulo each power of $p$.
Consider a term like $K^r A(1) J^{-s}$ and the
problem of making it integral at the prime $p$ for
sufficiently large $r$.  In order for it to be
integral, it suffices that it be so modulo the power
of $p$ dividing the determinant of $J^s$.  Modulo this
power of $p$, increase $r$ until we may replace
$K^r$ by $K_1$ and use
$$
K^r A(1)= U J_1 =U_1 J^s
$$
and
$$
K^r A(1)J^{-s} det(J)^{s} \equiv U_1 det(J)^{s}\pmod{det(J)^{s}}
$$
This guarantees the left hand side is divisible by
$det(J)^{-s}$, so that fractions in $$K^r A(1)J^{-s}$$
have no denominators $p$.  Taking all these primes
means we have no denominators at all.

A special case is the first equation
$$B(1)=J^{n(1)}A(1)^{-1}.$$
For this we allow a replacement similar to the above
of $A(1)$ by $K^rA(1)$ which won't affect solvability
of the other equations.
\end{proof}

 \begin{proposition}
Consider the matrix
\begin{eqnarray}
&&K=\left( \begin{array}{cccc}
0 &0 &1 &1\\
1 &0 &0 &0\\
0 &1 &0 &0\\
0 &0 &1 &0\end{array}\right)
\end{eqnarray}
Let $A(1)=K^{20}(K-1),J=KA(1)$ which is nonnegative.  Then
the matrices $J,K$ give isomorphic AF-algebras and
are unimodular with irreducible characteristic
polynomial,
but no powers of them are shift equivalent.
\end{proposition}

\begin{proof}
 These matrices come from \cite{KRW}, Ex.4.1
and \cite{KR92}, main theorem.  It can be checked they
are unimodular and are units in the field generated
by a root of the characteristic polynomial of
$K$.  By diagonalizing the field, it follows that
multiplication by $A(1)$ sends $V_1,w_1$ to $V_2,w_2$
in the notation of Lemma 3, so that the isomorphism
conditions are satisfied.  The eigenvalues of $J,K$ can
be identified with the matrices themselves and
their conjugates, under the map sending $K$ to its
maximal eigenvalue.  If powers of $J,K$ were shift
equivalent, then the
maximal eigenvalues would correspond to maximal
eigenvalues up to powers.  Hence we would have some
equation $K^r=(K^{21}(K-1))^s$.  But $K,K-1$ are independent
in the group of units of this degree 4 field, so there
can be no such equation.
\end{proof}

\begin{example}
\rm
Since the characteristic polynomial of a matrix of the form (18)
has the form
$$
 t^N - m_1t^{N-1} - m_2t^{N-2} -\ldots- m_{N-1}t - m_{N}
$$
it follows that two such matrices are conjugate (over $\mathbb{C}$)
if and only if they are equal. Proposition 5 shows
that two such matrices may satisfy the isomorphism
condition (14)--(15) without being equal. Let us mention an example
of two matrices satisfying (14)-(15) such that their twelfth powers
are conjugate over $\mathbb{C}$, but which are not shift equivalent:

\begin{equation}
J=\left( \begin{array}{ccccc}
 1 &1 &0 &0 &0 \\
 0 &0 &1 &0 &0 \\
 0 &0 &0 &1 &0 \\
 0 &0 &0 &0 &1 \\
 1 &0 &0 &0 &0\end{array} \right) ,\qquad
K=\left( \begin{array}{ccccc}
0 &1 &0 &0 &0 \\
0 &0 &1 &0 &0 \\
1 &0 &0 &1 &0 \\
1 &0 &0 &0 &1 \\
1 &0 &0 &0 &0\end{array} \right) .
\end{equation}

Their joint Perron-Frobenius eigenvalue is the real root of
 $t^3 - t -1$
and since they are unimodular, the dimension groups
are isomorphic by \cite[Corollary 6.2]{BJO98}. One computes that their
spectra are nondegenerate, and their twelfth powers
have the eigenvalue 1 with multiplicity 2 in addition
to the twelfth powers of the roots of the equation above.
Thus these powers are conjugate over $\mathbb{C}$.
But their 12th powers also seem to be not
shift equivalent.  The steps involved were to make the
two matrices block triangular using the 2 eigenvalues 1 and their
eigenvectors over the integers.  Then any conjugacy of the
matrices over the integers must also be triangular in this
form and for the $3\times 3$ block we must have a unit of the
field generated by a root of $z^3-29z^2-6z-1$, and the
units are generated by our matrix say $U_1$ and
$(7/55)U_1^2-(10/11)U_1-(57/55)$.   We have an equation for
the 3 nonzero blocks of the block triangular matrix
giving a shift equivalence, and modulo 7 one block in effect
cancels.  For the other we calculate up to scalar multiples
64 units reduced modulo 7 and check for each
that the equation is impossible.  We expect that what happens for
12th powers should also be true for all higher powers by a result in
\cite{KR79}.

Let us consider a more clearcut example:

\begin{equation}
J=\left( \begin{array}{cccccc}
 0 &1 &0 &0 &0 &0 \\
 1 &0 &1 &0 &0 &0 \\
 0 &0 &0 &1 &0 &0 \\
 0 &0 &0 &0 &1 &0 \\
 1 &0 &0 &0 &0 &1 \\
 1 &0 &0 &0 &0 &0\end{array} \right) ,\qquad
K=\left( \begin{array}{cccccc}
0 &1 &0 &0 &0 &0 \\
0 &0 &1 &0 &0 &0 \\
0 &0 &0 &1 &0 &0 \\
1 &0 &0 &0 &1 &0 \\
2 &0 &0 &0 &0 &1 \\
1 &0 &0 &0 &0 &0\end{array} \right) .
\end{equation}
The characteristic polynomials are
$$
\det(t\openone\! -\! J)=(t^3\!-t-1)(t^3\!+1) , \quad
\det(t\openone\! -\! K)=(t^3\!-t-1)(t^3\!+t+1) 
$$
and hence the Perron-Frobenius eigenvalue is the real root 
of $t^3 - t -1$
in both cases, and the associated dimension groups are isomorphic
by \cite[Corollary 6.2]{BJO98}.
All six roots of each of the two polynomials are distinct.
But we note that the three additional roots of the first polynomial
all have modulus one, and since all the roots of the latter 
polynomial has moduli different from one, $J^n$ is not conjugate
to $K^m$ over $\mathbb{C}$ for any positive powers $n, m$.

 We now give an example of two matrices $J,K$ of the form (18)
which are distinct, but have shift equivalent
second powers, so in particular the equivalence relation 4
is satisfied. To begin with let $J,K$ be the 6 by 6 matrices
$$
J=\left( \begin{array}{cccccc}
6 &1 &0 &0 &0 &0\\
16 &0 &1 &0 &0 &0\\
197 &0 &0 &1 &0 &0\\
90 &0 &0 &0 &1 &0\\
2200 &0 &0 &0 &0 &1\\
12000 &0 &0 &0 &0 &0\end{array} \right) , \qquad
K=\left( \begin{array}{cccccc}
8 &1 &0 &0 &0 &0\\
2 &0 &1 &0 &0 &0\\
97 &0 &0 &1 &0 &0\\
370 &0 &0 &0 &1 &0\\
3400 &0 &0 &0 &0 &1\\
12000 &0 &0 &0 &0 &0\end{array} \right) ,
$$
respectively. Their characteristic polynomials are
$$
\det(t\openone\! -\! J)=(t-10)(t+3)(t^2-4t+16)(t^2+5t+25),
$$
$$
\det(t\openone\! -\! K)=(t-10)(t+3)(t^2+4t+16)(t^2-5t+25),
$$
respectively. The squares $J^2,K^2$ of these matrices
have the same characteristic
polynomial
$$
(t-100)(t-9)(t^2+16t+256)(t^2+25t+625)
$$
which has six distinct roots. Thus there are nontrivial
intertwiners between each of $J,K$ and each of the matrices
$$
\left( \begin{array}{c}
 100\end{array} \right) , \qquad
\left( \begin{array}{c}
 9\end{array} \right) , \qquad
\left( \begin{array}{cc}
 -8 &14\\
 -14 &-8\end{array} \right) , \qquad
\left( \begin{array}{cc}
 -12 &67\\
 -7 &-13\end{array} \right) , 
$$
 which are all irreducible over $\mathbb{Q}$. At the
outset, these intertwiners may be taken to have
matrix elements in a finite field extension of order
1 or 2 over $\mathbb{Q}$, but then by linearity they
may be taken to have matrix elements in $\mathbb{Q}$.
Thus $J^2$ and $K^2$ are  similar over $\mathbb{Q}$, and
hence they are elementary shift equivalent over 
$\mathbb{Q}$.  It follows by an argument using some rational
shift equivalences with denominator $d$, that
if we replace $J,K$ by
$$
J=\left( \begin{array}{cccccc}
6d &1 &0 &0 &0 &0\\
16d^2 &0 &1 &0 &0 &0\\
197d^3 &0 &0 &1 &0 &0\\
90d^4 &0 &0 &0 &1 &0\\
2200d^5 &0 &0 &0 &0 &1\\
12000d^6 &0 &0 &0 &0 &0\end{array} \right) , \qquad
K=\left( \begin{array}{cccccc}
8d &1 &0 &0 &0 &0\\
2d^2 &0 &1 &0 &0 &0\\
97d^3 &0 &0 &1 &0 &0\\
370d^4 &0 &0 &0 &1 &0\\
3400d^5 &0 &0 &0 &0 &1\\
12000d^6 &0 &0 &0 &0 &0\end{array} \right) ,
$$
for an appropriate $d$ in $\mathbb{N}$
we obtain shift equivalence of the squares $J^2,K^2$ 
over the integers. (The details of the argument are written
out in e.g. \cite[Section7]{BJO98}.) But then $J^2,K^2$ are shift
equivalent over the non-negative integers by
\cite[Theorem 2.1]{KR79}, \cite[Theorem 7.3.6]{LM95}.
(As explained after (18), $J,K$
themselves cannot be shift equivalent.)

Let us end with an example which is maximally clearcut,
and clears the way for Theorem~10:
\begin{equation}
J=\left[6 \right] ,\qquad
K=\left[12 \right] .
\end{equation}
We leave it to the reader to decide that $J,K$ are $C^*$-equivalent,
i.e. (15) holds, but no powers are shift equivalent or even conjugate.

\end{example}

We next give a necessary condition for $C^*$-equivalence.
Recall that if $\lambda$ is an algebraic number, then
$\mathbb{Q}[\lambda]$ is a field. ( $\mathbb{Q}[\lambda]$ 
is by definition
a ring, all polynomials in $\lambda$ over $\mathbb{Q}$,
and multiplicative inverses of 
nonzero elements 
exist because of the Euclidean algorithm in $Q[t]$: Let $f(t)$ be 
the prime polynomial corresponding to $\lambda$ (so
that $\mathbb{Q}[\lambda]=\mathbb{Q}[t]/f[t]$) and
let $g(t)$ be any polynomial nonzero modulo $f$.  Then there
exist polynomials $r,s\in Q[t]$ such that $rg+sf=1$, hence
modulo $f$, $r$ is an inverse of $g$.) The term "prime" in
this field still means a prime in the associated
subring of algebraic integers.

\begin{theorem}
Suppose that the two non-negative square 
$d\times d$ matrices $J,K$ are nonsingular and primitive.  Let
$\pi_1,\pi_2$ the projections from all row vectors
 to the maximal eigenspaces of $J,K$.
Then a necessary condition
for isomorphism of the AF algebras determined by $J,K$
 is that 
\vskip0.3cm
(1) (\cite[Corollary 3.5]{BD95}) the fields $\mathbb{Q}[\lambda_1]$, 
$\mathbb{Q}[\lambda_2]$ generated by
the maximal eigenvalues $\lambda_1,\lambda_2 $
of $J,K$ are the same and $\lambda_1,\lambda_2$ are products
of the same primes over this field,
 \vskip0.3cm
 (2) the dimension group
quotients denoted $\pi_i(G)=G^{Per}$ in \cite{BMT87} (quotients by
nonmaximal eigenspaces) are isomorphic as modules over
 $\mathbb{Z}[1/\lambda_i^n]$ for sufficiently large $n$ .
 \vskip0.3cm
This is
a sufficient condition if the characteristic polynomials
are irreducible, otherwise Theorem 7 gives
some additional necessary conditions.
\end{theorem}

\begin{proof}
  
First we argue for the necessity. Let the maximal and
nonmaximal eigenspaces of $J,K$ be $<w_1>,V_1,<w_2>,V_2$
as in Lemma 3. The matrix $A(1)$
must map $V_1$ to $V_2$ nontrivially, and send $w_1$ to
$w_2$ nontrivially by Theorem 6.  Since $A(1)$ is rational,
it commutes with Galois actions among the different 
conjugates of the maximal eigenvalue.
This already implies that the two fields must be the same (they
have the same set of nontrivial Galois actions under some
finite Galois extension containing both).
It follows essentially by Theorem
7 that the prime factors of the two maximal eigenvalues
must be the same over algebraic number fields.
The intersection of all conjugates of $V_1$ goes to the 
corresponding intersection for $V_2$ by this isomorphism
consistent with Galois action, so that the quotient of the
rational dimension group by nonconjugate eigenvalues is
mapped from the one to the other.  The mapping $A(1)$ taken
over the field $\mathbb{Q}[\lambda]$ gives an isomorphism between
the maximal eigenspaces.  This will be multiplication by some element
 $x\in \mathbb{Q}[\lambda]$. 
 Let $G_1,G_2$ be
$\pi_1(\mathbb{Z}^d),\pi_2(\mathbb{Z}^d)$.
Multiplication by $x$ will take $G_1$ into $G_2$, by the
effect of $A(1)$, and its inverse will do the reverse,
up to multiplication by powers of the primes in $\lambda_i$
which represent multiplication by $J,K$.  Therefore 
\vskip0.3cm
(DG) $xG_1\subset G_2, \lambda_1^{e_0}x^{-1}G_2\subset G_1$
for some positive $e_0$.
 \vskip0.3cm
Recall that the dimension groups of $J,K$ can be viewed as the
direct limit of $\mathbb{Z}^d$ sent to itself by $J,K$; if we embed them
in the maximum eigenspace this direct limit is equivalent to
making $\lambda_1$ or $\lambda_2$, or all the primes in them
invertible so that (DG) implies isomorphism.
This implies that the dimension groups of the matrices $J,K$
are isomorphic as modules over $\mathbb{Z}[1/\lambda_i^n]$ for 
sufficiently
large $n$. (This condition does not depend on $i$.)
This proves necessity.

Now assume the characteristic polynomials of $J,K$ are irreducible
and (1),(2) hold, hence (DG).  Choose such an isomorphism as in
(DG), and
adjust its sign so that on the maximal eigenvector it is
positive.  Expanding out the coefficients gives a
map $A(1)$ over the integers which preserves all the
conjugates of the rational eigenspace.  The effect of $A(1)$
on all vectors over $\mathbb{Q}$ is isomorphically mirrored
in its effect on vectors over $\mathbb{Q}[\lambda]$ on the maximal
eigenspace there.  This fact together with Lemma 3 ensures
positivity of the other $A(i),B(i)$.  To get positivity
of $A(1)$ we replace it by some product
$K^{r_1}A(1)J^{r_2}$ and note that its effect 
on the maximum row eigenvector
is positive and this dominates the product
asymptotically. The assumption on the 
primes and (DG), implies that for all
$n$ there exists an $m$ such that $K^m A(1) J^{-n},K^m A(1)^{-1}J^n$
involves only nonnegative powers of all the primes and is divisible
by any given power of each $\lambda_i$.
This implies that all the $A(i),B(i)$ exist over the integers.

\end{proof}

\begin{acknowledgements}
We are indepted to Akitaka Kishimoto for pointing
out that the first problem raised in this paper really is a
problem, to Danrun Huang for establishing the 
contact between the two teams of authors, and to the
referee for making us aware of \cite{SV98}.
\end{acknowledgements}

\end{document}